\documentclass[12pt]{article}
\begin{document}
\begin{center}
\large{\bf{A generalization of Kakutani's splitting procedure}}
\smallskip

{\large{\bf{Aljo\v{s}a Vol\v{c}i\v{c}}}}
\end{center}
\bigskip

\begin{abstract}
In this paper we will study 
uniformly distributed sequences of partitions of $[0,1]$, a concept
which has been introduced by Kakutani in 1976. We will construct new
families of u.d. sequences of partitions and study relations with the
classical concept of uniformly distributed sequences of points.

AMS classification: 11K06, 11K45, 28D05, 60B10, 60F15, 60F20

\end{abstract}

\section {\bf Introduction} 
\smallskip

If we split the interval $[0,1]$ in two equal parts, and then in four equal parts and keep going, splitting it in $2^n$ equal parts at the $n$-th step, we get the sequence of dyadic (or binary) partitions which are finer and finer. The procedure is, in an obvious sense, uniform. 

If on the other hand we take any $\alpha  \in ]0,1[$ different from one half and split $[0,1]$ in $[0,\alpha]$ and $[\alpha,1]$, and follow a similar procedure as before, splitting each interval of the previous partition in two parts, left and right, proportional to $\alpha$ and $\beta=1-\alpha$ respectively, then we get the so-called $\alpha$-dyadic partitions, whose behavior depends of course on $\alpha$. An obvious observation is that the points of the partitions are now not evenly distributed, since $[0,\alpha]$ and $[\alpha,1]$ contain the same number of points of any subsequent $\alpha$-dyadic partition, in spite of the fact that they do not have the same measure.

In [15] Kakutani introduced  the following interesting variant of the splitting procedure just described. 
\bigskip

\noindent {\bf Definition 1.1} If $\alpha \in ]0,1[$ and
$\pi=\{[t_{i-1} , t_i] : 1\le i \le k\}$ is any partition of $[0,1]$, its {\it $\alpha$-refinement} (which will be denoted by $\alpha \pi$) is obtained by
subdividing  {\it only the interval(s) of $\pi$ having maximal length} in two parts, proportional to $\alpha$ and $\beta=1-\alpha$, respectively. 
\bigskip

Kakutani's sequence of partitions \{$\kappa_n$\} is
obtained by successive $\alpha$-refine\-ments of the trivial partition
$\omega =\{[0,1]\}$. 

\noindent For example, if $\alpha <\beta$, $\kappa_1=
\{[0,\alpha] ,[\alpha, 1]\}$, $\kappa_2=\{[0,\alpha] ,[\alpha, \alpha +\alpha \beta],[\alpha +\alpha \beta, 1]\}$, and
so on. 

Note that for $\alpha= \frac{1}{2}$ we get just the sequence of dyadic partitions.
\smallskip

To express things with the appropriate precision, we need the following definition.
\smallskip

\noindent {\bf Definition 1.2} Given a sequence of
partitions $\{\pi_n\}$ of $[0,1]$, with 
$$\pi_n=\{[t_{i-1}^n , t_i^n] :1\le i \le k(n)\} \,,$$ 
we say that it is {\it uniformly distributed (u.d.)}, if for
any continuous function $f$ on $[0,1]$ we have 
$$\lim_{n\rightarrow \infty}\frac{1}{k(n)} \sum_{i=1}^{k(n)} f(t_i^n)= \int_0^1 f(t)\,d
t\,.
 \eqno (1)$$ 
\smallskip

A simple but useful observation is that uniform distribution of a sequence of partitions $\{\pi_n\}$ is equivalent to the weak convergence to $\lambda$ (the Lebesgue measure on $[0,1]$) of the associated sequence of measures $\{\nu_n\}$, with
$$\nu_n=\frac{1}{k(n)} \sum_{i=1}^{k(n)} \delta_{t_i^n}\,,$$
where $\delta_t$ denotes the Dirac measure concentrated in $t$.
\smallskip

This reformulation allows to use classical compactness criteria for weak convergence (see [19]).

We will use in the sequel another standard argument: a sequence of partitions is u.d.  if and only if (1) holds for a family of Riemann integrable functions $\cal{G}$ such that its linear span is dense in the class of the continuous functions on $[0,1]$ with respect to the $\|\cdot \|_{\infty}$ norm ([16], Chapter 1).
\smallskip

We can state now Kakutani's main result.
\bigskip

\noindent {\bf Theorem 1.3} {\it For any $\alpha \in ]0,1[$ the
sequence of partitions \{$\kappa_n$\} is uniformly distributed.}
\bigskip

This fact got a considerable attention in the late seventies and early eighties, when
other authors provided different proofs of Kakutani's theorem  [1] and also
proved several stochastic versions, in which the intervals of maximal length are split according to certain probability distributions ([17], [18], [25], [22], [23] and [4]). The expository work [2] gives a complete overview on the results obtained before 1984.
\smallskip

The paper [8] extends the notion of uniform distribution of a sequence of partitions to probability measures on complete separable metric spaces.

More recently some new questions and ideas revived the interest for the subject and few items have been added to the pertinent bibliography ([20], [6], [7]). The second of these papers in particular extends Kakutani's splitting procedure to higher dimension with a construction which is intrinsically higher-dimensional. Further results are on the way, like for instance the extension of the theory to certain fractals and the explicit computation of the discrepancy of a particular class of sequences of partitions.

This paper will contribute further results to this subject.
\smallskip

In Section 2 we generalize Kakutani's procedure by splitting the longest interval into a finite number of parts homothetically to a given partition $\rho$ and prove that $\{\rho^n\omega\}$ is uniformly distributed.

In Section 3 we associate to a given u.d. sequence of partitions random reorderings of the points determining the partitions, proving that almost surely they provide uniformly distributed sequences of points. 
This answers a question posed by David Fremlin.
\smallskip

Section 3 provides an explicit  link between our area of investigation and the theory of uniformly distributed sequences of points, initiated by Hermann Weyl in [26]. We believe our results will have applications in quasi-Monte Carlo methods, because new simple u.d. sequences of low discrepancy (to be discussed in another paper) can be generated by our methods, and combined with the approach initiated by van der Corput (in [24] and seven subsequent papers numbered from I to VIII), Halton [10], [11] and Hammersley [12].
\smallskip

An excellent and by now classic exposition of results on uniformly distributed sequences of points is [16]. For more recent results and applications  see for instance  [9], in particular Chapter 3.
\bigskip

\noindent {\it Acknowledgements} The autor wants to thank David Fremlin and Pietro Rigo for helpfull discussions.
\smallskip

\section {\bf $\rho$-refinements}
\smallskip

In this section we will geneneralize Kakutani's splitting procedure and produce a new class of u.d. sequences of partitions.
\smallskip

Let us give the following definition.
\bigskip

\noindent {\bf Definition 2.1} Consider any non trivial finite partition $\rho$ of $[0,1]$. We will keep it fixed for the whole section. The {\it $\rho$-refinement} of a partition
$\pi$ of $[0,1]$ (which will be denoted by $\rho \pi$) is obtained by
subdividing all the intervals of
$\pi$ having maximal length positively (or directly) homothetically to $\rho$.
\smallskip

Obviously, if $\rho =\{[0, \alpha],[\alpha, 1]\}$, then the $\rho$-refinement is just
Kakutani's $\alpha$-refinement. 

As in Kakutani's case, we can iterate the splitting procedure. The $\rho$-refinement of $\rho \pi$ will be denoted by $\rho^2 \pi$, and the meaning of $\rho^n \pi$, for $n\in I\!\!N$, does not need any explanation.

In this section we will  prove that the sequence $\{\rho^n\omega\}$ of successive $\rho$-refinements of the trivial partition $\omega$ is u.d..
\bigskip

\noindent {\bf Remark 2.2} It is worth noticing that it is necessary to put some restriction on the partition  $\pi$ (even in the simplest case of the Kakutani splitting procedure)  if we hope for uniform distribution of $\{\rho^n\pi\}$.

It would be interesting to find significant sufficient conditions on $\pi$ in order to obtain the uniform distribution of $\{\rho^n\pi\}$ even for the case of Kakutani's splitting procedure.
\smallskip

The following simple example shows where lies the problem.
\smallskip

Let $\pi=\{[0,\frac{2}{5}], [\frac{2}{5},1]\}$ and consider $\rho=\{[0,\frac{1}{2}], [\frac{1}{2},1]\}$. It is clear that the $\rho$-refinement operates alternatively on $ [\frac{2}{5},1]$ and $[0,\frac{2}{5}]$,  and that the subsequences $\{\rho^{2n}\pi\}$ and  $\{\rho^{2n+1}\pi\}$ converge to measures which attribute to $[0,\frac{2}{5}] $ the values $\frac{5}{4}$ and $\frac{5}{6}$, respectively, so $\{\rho^n \pi \}$ does not converge.
\smallskip

Let us now fix some notations and recall some preliminarily facts which will be used to prove Theorem 2.7.
\smallskip

Let $\rho = \{[r_{i-1},  r_i]:  1\le i \le k \}$ be the fixed partition 
of $[0,1]$. Denote by $\alpha_i =r_i-r_{i-1}$, for $1\le i \le k$, the
lengths of the $k$ intervals of $\rho$.

Let us denote by $[\rho]^n$ the  so-called $n$-th $\rho$-adic partition of $[0,1]$, obtained from 
$[\rho]^{n-1}$ (where $[\rho]^1=\rho$) by subdividing {\it all} its $k^{n-1}$ 
intervals positively homothetically to $\rho$. If an interval belongs to $[\rho]^n$, we will say that it has {\it rank} $n$.

The $k$ intervals $[\sum_{h=1}^{i-1}\alpha_h,\sum_{h=1}^{i}\alpha_h]$,
$1\le i \le k$,  of $\rho$ will be denoted by $I(\alpha_i)$. 
If $[y_{j-1},y_j]=I(\alpha_{i_1}\alpha_{i_2}\dots\alpha_{i_{n-1}})$,
$1\le j \le k^{n-1}$, is a generic interval of rank $n-1$,
the $k^n$ intervals of $[\rho]^n$ are
$$I=I(\alpha_{i_1}\alpha_{i_2}\dots\alpha_{i_{n-1}}\alpha_{i_n})$$
$$=[y_{j-1}+(y_{j}-y_{j-1})\sum_{h=1}^{i_n-1}\alpha_h,\,
y_{j-1}+(y_{j}-y_{j-1})\sum_{h=1}^{i_n}\alpha_h]\,.$$

Here the $\alpha_i$'s are understood as symbols (following [1]), not as numbers. So if for instance $\alpha_1=\alpha_2=\alpha_3=\frac{1}{3}$, $I(\alpha_1\alpha_2)=[\frac{2}{9}, \frac{1}{3}]$, while $I(\alpha_2\alpha_1)=[\frac{1}{3},\frac{4}{9}]$.

Note that $\lambda(I(\alpha_{i_1}\alpha_{i_2}\dots\alpha_{i_{n-1}}\alpha_{i_n}))=
\prod_{m=1}^n\alpha_{i_m}$.
\smallskip

Let $X=\{ \alpha_1, \alpha_2, \dots, \alpha_k \}$ and let 
$\sigma$ be the probability on $X$ such that $\sigma (\{\alpha_i\})=\alpha_{i}$ for $1\le i\le k$.

Put $X_m=X$ and $\sigma_m =\sigma$ for any $m\in I\!\!N$.

Denote by $Y$ the countable product $\prod_{m=1}^{\infty}X_m$. On $Y$ we
consider the usual product probability $\mu$.  If
$C=C(\alpha_{i_1}\alpha_{i_2}\dots\alpha_{i_n})=\prod_{m=1}^{\infty}X_m'$,
where $X_m'=\{\alpha_{i_m}\}$ for $m\le n$ and $X_m'=X$ for $m>n$, is a cylinder set, then 
$\mu(C)=\prod_{m=1}^n\alpha_{i_m}$.

To every point $t\in [0,1]$ we can associate a sequence
$\{\alpha_{i_m}\}$ such that
$$t\in
\cap_{m=1}^{\infty}I(\alpha_{i_1}\alpha_{i_2}
\dots\alpha_{i_{m-1}}\alpha_{i_m})\,.$$

There is an (expected) ambiguity to be taken care of: there are two such sequences $\{\alpha_{i_m}\}$  for the countable set of points belonging to the endpoints of some interval of $ [\rho]^n$; in that case
we associate to $t$ the sequence for which definitely $\alpha_{i_h}=\alpha_1$.
\bigskip

This defines a $1-1$ mapping between $[0,1]$ and a subset $Y'$ of $Y$, obtained by removing from Y the countable set of sequences $\{ \alpha_{i_m}\}$ such that definitely $\alpha_{i_m}=\alpha_k$. Note that $\mu (Y\setminus Y')=0$.

This mapping is measure preserving if we take on $[0,1]$  the Lebesgue
measure and on $Y'$ the restriction of $\mu$. This follows immediately observing that cylinder sets
$C$ and $\rho$-adic intervals $I$ bearing the same indices have the same measure.
\smallskip

If $I$ and $J$ are two disjoint subintervals of $[0,1]$ having the same length, and $J=I+c$, say, the function $f_{I,J}$ which takes values $x+c$ on $I$, $x-c$ on $J$ and $x$ otherwise, is measure preserving.

Let us denote by $\cal{F}$ the family of all such functions  $f_{I,J}$, where now $I$ and $J$ are disjoint $\rho$-adic intervals of the same length. The intervals $I$ and $J$ need not have the same rank. The next proposition proves an important property of $\cal{F}$.
\smallskip

\noindent {\bf Lemma 2.3} {\it The family $\cal{F}$ is ergodic.}
\smallskip

\noindent {\bf Proof.} We will exploit the measure isomorphism described above. Denote by $\cal{F'}$ the corresponding family of transformations on $Y'$. When $f=f_{I,J}$ and $I$ and $J$ do have the same rank, the corresponding function  $f'$ on $Y'$ is a permutation of a finite number of coordinates and it  preserves the product measure $\mu$. The family of these transformations is ergodic by a theorem due to Hewitt and Savage ([14]). Since it is contained in $\cal{F'}$, the latter is also ergodic. By the isomorphism between $(Y',\mu)$ and  $([0,1], \lambda)$, we conclude that $\cal{F}$ is ergodic, too.
\bigskip

For the partition $\rho^n \omega$, let $A_n$ denote the length of the longest interval and $a_n$ the length of the shortest interval.  We have the following estimate.
\bigskip

\noindent {\bf Lemma 2.4} {\it For any $n$, 
$$a_1A_n \le a_n\,.$$}

\noindent {\bf Proof.}  Since $A_n<1$, the inequality is true for $n=1$. Proceed now by induction. Suppose 
$a_1A_{n-1} \le a_{n-1}$. There are two possibilities: either $a_n= a_{n-1}$ or $a_n< a_{n-1}$. In the first case, since $A_n<A_{n-1}$ for any $n$, 
$$a_1A_n<a_1A_{n-1} \le a_{n-1}=a_n\,.$$
In the second case the shortest interval of $\rho^n\omega$ results from the splitting of the longest one of $\rho^{n-1}\omega$, in other words $a_n=a_1A_{n-1}$, so the conclusion follows in either case.
\bigskip

\noindent {\bf Lemma 2.5} {\it If $\pi_n= \rho^n\omega$, then}
$$\lim_{n\rightarrow \infty}  \mbox{diam} \,\, \pi_n=0\,.$$

\noindent {\bf Proof.}  With the notation introduced above, obvioulsy $a_n\le \frac{1}{n}$ for any $n$. Apply now the previous lemma to get
$$\mbox{diam} \,\,\pi_n=A_n\le \frac{a_n}{a_1}\le \frac{1}{n a_1}\,.$$
\bigskip

\noindent {\bf Remark 2.6}  An immediate consequence of the previous lemma is that if we consider the family of the characteristic functions (which are Riemann integrable) of all the intervals belonging to the partitions $\rho^n \omega$, $n\in I\!\!N$,  its linear span is dense in the class of all continuous functions on $[0,1]$ with respect to the $\| \cdot \|_{\infty}$ norm. 
\smallskip

We are now ready for proving the main result of this section.

\bigskip

\noindent {\bf Theorem 2.7} {\it The sequence $\{\rho^n \omega\}$ is uniformly distributed.}
\smallskip

\noindent {\bf Proof.} It is well known that the set of all Borel probability measures on $[0,1]$, with the topology associated to the weak convergence, is  metrizable and compact ([19], Theorem 6.4). Let us denote by $\{\nu_n\}$ the sequence of measures associated to the partitions $\{\rho^n \omega\}$ as defined in Section 1. Then $\{\nu_n\}$ admits weakly convergent subsequences. All we have to prove is that any such subsequence converges to $\lambda$.

The plan of the proof is the following: let $\nu$ be the limit of a convergent subsequence of $\{\nu_n\}$. Then Lemmas 2.4 and 2.5 provide lower and upper bounds for the Radon-Nikodym derivative $\frac{d\nu}{d\lambda}$. The ergodicity of the family $\cal F$ tells us finally that this derivative has to be constant  almost everywhere and therefore $\nu=\lambda$.
\smallskip

Every interval $J$ belonging to some $\rho^n \omega$ belongs to some $[\rho]^m \omega$, but it is also true, viceversa, that every interval $I\in [\rho]^m \omega$ sooner or later appears in one of the partitions $\rho^n \omega$. This is due to the fact that 
$$s= \sup \{r: J_r\in \rho^r\omega, I\subset J_r \}$$
is well defined, since $\rho$-adic intervals are either disjoint or one contains the other.
Moreover the diameter of $\rho^n \omega$  tends to zero. This implies that $I=J_s$.

Let now $J$ be any $\rho$-adic interval and suppose that $m\in I\!\!N$ is such that every $\rho^n \omega$, for $n\ge m$, subdivides $J$. 

If $J$ is subdivided by $\rho^n \omega$ into $k$ intervals, then obviously
$$ka_n \le \lambda(J) \le k A_n\,,$$
where $a_n$ and $A_n$ are the quantities considered in Lemma 2.4.

If $k(n)$ is the number of intervals of $\rho^n \omega$,  the previous double inequality can be rewritten, in terms of the measure $\nu_n$ (by Lemma 2.4)  as
$$\frac{\lambda (J)}{k(n)A_n}\le \nu_n(J)\le \frac{\lambda (J)}{k(n)a_n}\,,$$
for any $n\ge m$. It follows  that, for all $n$ sufficiently large,
$$a_1\lambda (J) \le \frac{a_1 \lambda (J)}{k(n)a_n}\le \nu_n (J) \le \frac{\lambda (J)}{k(n)a_1 A_n}\le \frac{\lambda (J)}{a_1}\,.$$

Suppose now that a subsequence $\{\nu_{n_k}\}$ converges weakly to $\nu$. From the last inequality we can conclude that for any $\rho$-adic interval $J$ it is
$$a_1\lambda (J) \le  \nu (J)  \le \frac{\lambda (J)}{a_1}\,.$$
Hence, by Remark 2.6, the same inequality holds for any Borel set $B$.

Therefore $\lambda <\!\!< \nu <\!\!< \lambda$ and if we denote by $\frac{d\nu}{d\lambda}$ the Radon-Nikodym derivative of $\nu$ with respect to $\lambda$, we have $\lambda$-a.e.
$$a_1\le \frac{d\nu}{d\lambda} \le \frac{1}{a_1}\,.$$

Observe that if $I$ and $J$ are two intervals which have the same length and belong to some $\rho^n\omega$ (they need not have the same rank), then our splitting procedure behaves on them in the same way. This implies that $\nu(I)=\nu(J)$.

Therefore $\nu$ is invariant with respect to the family of functions $\cal{F}$ considered in Lemma 2.3, where we proved that it is ergodic. It follows that $\frac{d\nu}{d\lambda}$ is constant almost everywhere, i.e. $\nu=\lambda$. In other words, the sequence $\{\rho^n\omega\}$ is u.d..

\section{\bf Associated uniformly distributed sequences of points}

In this section we will study the relation between u.d. sequences of
partitions and u.d. sequences of points.

Let us recall the following definition. 
\bigskip

\noindent {\bf Definition 3.1} A sequence of points $\{x_n\}$ in $[0,1]$ is said to be {\it uniformly distributed (u.d.)}зк, if 
$$\lim_{N\rightarrow \infty} \frac{1}{N}\sum_{i=1}^Nf(x_i)=\int_0^1 f(x)\,dx\,, $$
for every continuous function on $[0,1]$.
\smallskip

Suppose $\{\pi_n\}$ is a u.d. sequence of partitions, with $\pi_n = \{[t_{i-1}^n , t_i^n]: 1\le i \le k(n)\}$. Then it is natural to ask whether the points $\{t_i^n\}$ determining the partition $\pi_n$, for $1\le i \le k(n)$ and $n\in I\!\!N$, can be rearranged in order to get a u.d. sequence of points.

There are of course many ways of reordering the points $\{t_i^n\}$. A natural restriction to make  is that we first reorder all the points defining $\pi_1$, then those defining $\pi_2$, and so on.  We will limit ourselves to such kind of reorderings  and call them {\it sequential}. But still there is a lot of freedom left  in rearranging the points $\{t_i^n\}$, for $n$ fixed, and not all of them produce the desired result.
\smallskip

Consider again the dyadic sequence of partitions. Put
$I_h^n=[\frac{h-1}{2^n}, \frac{h}{2^n}]$, for $n\in I\!\! N$ and $h=
1, 2, \dots , 2^n$. Let us order, for any fixed $n$, the left end-points of the $I_h^n$'s by magnitude. We obtain in this way the so-called {\it lexicographic} order which does not give a u.d. sequence, as it is easy to check. 

It is known from the literature that the (by far) best thing one can do is the following (see  [16], Theorem 3.5).
Any positive integer $k$ can be represented uniquely by its dyadic expansion
$k=\sum_{i=1}^s a_i2^i$, with $a_i\in \{0,1\}$, $a_s=1$.

Set now, for any $k\ge 1$,  $x_k=\sum_{i=1}^s a_i2^{-i-1}$. This is called the van der Corput's sequence.
It is a well known result that $\{x_k\}$ is u.d. and that it has optimal discrepancy (for more information, see also [9], Chapter 1 and also the brief discussion in the last section of this paper).
\smallskip

So the lexicographic order is no good, van der Corput's method is excellent. What can be said about the generic sequential reordering? We will give here a probabilistic answer which applies to any u.d. sequence of partitions.
\smallskip

First we need a version of the strong law of large numbers for negatively correlated random variables. This result is attributed in [4] to Aleksander Rajchman and can be proved, with small modifications, along the lines of [4], Theorem 5.1.2.
\bigskip

\noindent {\bf Lemma 3.2} {\it Let $\{\varphi_n\}$ be a sequence of real, negatively correlated random variables with variances uniformly bounded by $V$ on the probability space $(W,P)$. 
Suppose moreover that
$$\lim_{i \rightarrow \infty} E(\varphi_i)=M\,.$$
Then
$$\lim \frac{1}{n}\sum_{i=1}^n \varphi_i=M \,\,\mbox{almost surely}\,.$$}
\bigskip

We will also need few simple properties of a particular family of discrete random variables. Let $\varphi$ be the random variable taking  with probability $\frac{1}{k}$ values in the sample space $W=\{w_i \in [0,1], 1\le i\le k\}$,  $k\ge 2$. We may, and do, assume that $w_{i-1}<w_i$ for $1\le i \le k$. Denote by $\varphi_i$ the value of $\varphi$ in the $i$-th draw from $W$, without replacement. Fix now $c\in ]0,1[$ and let $\psi_i=\chi _{[0,c[}(\varphi_i)$, where $\chi _{[0,c[}(\cdot)$ denotes the characteristic function of $[0,c[$. We have then the following simple property.
\bigskip

\noindent {\bf Proposition 3.3 } {\it The variances of the random variables $\psi_i$, $1\le i\le k$, are bounded by $\frac{1}{4}$. Moreover, the $\psi_i$'s are negatively correlated.}

\bigskip

Let us now prepare the setting for the main result of this section.
\smallskip

If $\{\pi_n\}$ is a sequence of uniformly distributed partitions of $[0,1]$, where $\pi_n=\{[t_{i-1}^n, t_i^n]: 1\le i \le k(n)\}$, the {\it sequential random} reordering of the points $t^n_i$ is a sequence (denoted by $\{\varphi_m\}$) made up by consecutive blocks of $k(1), k(2), \dots ,$ $ k(n), \dots $  random variables. The random variables of the $n$-th block have the same law and represent the drawing, without replacement, from the sample space $W_n=\{t_1^n, t_2^n, \dots , t_{k(n)}^n\}$, where each singleton has probability $\frac{1}{k(n)}$. 

Denote by $T_n$ the set of all permutations on $W_n$, endowed with the natural probability compatible with the uniform probability on $W_n$: each permutation $\tau_n \in T_n$ has probability $\frac{1}{k(n)!}$. 

Any sequential random reordering of the sequence $\{\pi_n\}$ corresponds to a random selection of $\tau_n\in T_n$, for each $n\in I\!\!N$. The permutation $\tau_n$ identifies the reordered $k(n)$-tuple of random variables $\varphi_i$ for $K(n-1) < i \le K(n)$, where  $K(n)=\sum_{i=1}^{n}k(i)$. The set of all sequential reorderings can be given therefore the natural product probability on 
$T=\Pi_{n=1}^{\infty}T_n$. 
\bigskip

\noindent {\bf Theorem 3.4 } {\it If $\{\pi_n\}$ is a uniformly distributed sequence of partitions of $[0,1]$, then the sequential random reordering of the points $t^n_i$ defining them is almost surely a uniformly distributed sequence of points in $[0,1]$.}

\noindent {\bf Proof. } Observe first that if $0<c<1$ and $\varphi_m$ belongs to the $n$-th block of $k(n)$ random variables, then
$$E(\chi_{[o,c[} (\varphi_m))=\frac{1}{k(n)}\sum_{i=1}^{k(n)}\chi_{[o,c[} (t_i^n)\,,$$
and this quantity tends to $c$, when $m$ and hence $n$ tend to infinity, since $\{\pi_n\}$ is u.d..

Observe that Proposition 3.3 applies to any set of $k(n)$ random variables $\psi_m=\chi_{[o,c[} (\varphi_m)$, for 
$K(n-1) < m \le K(n)$. It follows that the sequence $\{\psi_m\}$ is negatively correlated, since the correlation is negative if the random variables belong to the same block and is zero if they belong to different blocks, because they are independent. 

Let $\{c_h, h\in I\!\!N\}$ be a dense subset of $[0,1]$.  Fix $h\in I\!\!N$ and consider  the sequence 
$\{\chi_{[o,c_h [} (\varphi_m)\}$. We may apply to it Proposition 3.2 and get
$$\lim_{n \rightarrow \infty} \frac{1}{n}\sum_{i=1}^n \chi_{[o,c_h [} (\varphi_i)=c_h \,\,\,\,\mbox{a. s.}\, $$
for any $c_h$.

But this is a well known sufficient condition for uniform distribution (compare Exercises 1.1 and 1.3 in [15]).

\section{\bf Final remarks}

The two main results of this paper, Theorem 2.6 and Theorem 3.4, rise several questions.

The first result calls for explicit computation of the discrepancy of the sequence $\{\rho^n\omega\}$. This has been done for a class of these sequences and will appear elsewhere.

For a finite set of points $W=\{t_i\in [0,1] : 1\le i \le k\}$, the discrepancy is defined by
$$D_W=\sup_{0\le a<b\le 1} \left|\frac{\sum_{i=1}^k \chi_{[a,b[}(t_i)}{k}  - (b-a)\right|\,.$$

A natural problem is to estimate the behavior of the sequence $\{D_{W_n}\}$ (where $W_n$ is the set defining the partitions $\rho^n\omega$) when n tends to infinity.  It is easy to see that a sequence of partitions is u.d. if and only if $\{D_{W_n}\}$ tends to zero, when $n$ tends to infinity, but it is of particular interest to find partitions $\rho$ that such that the speed of convergence is as high as possible. Some classes have been found and will be a subject of a separate paper.

It is easy to see that if $W_n=\{0, \frac{1}{n}, \dots , \frac{n-1}{n}\}$, then $\{D_{W_n}\}=\frac{1}{n}$ and that this speed cannot be improved ([16], Chapter 2, Theorem 1.2 and Example 1.1).

On the other hand a classical result says that for u.d. sequences of {\it points}, the discrepancy behaves, at best, as $\frac{\log n}{n}$. This has been a long-standing open problem solved by W. M. Schmidt [21]. Sequences of points whose discrepancy is of this order are said to have {\it low discrepancy}. 
\smallskip

There are also several questions connected to Theorem 3.4.

It would be interesting to compare the discrepancy of $\{\pi_n\}$ and that of ``most" of its sequential random reorderings. In particular, if the discrepancy of the sequence of partitions behaves like $\frac{1}{n}$, what is the chance to get a low discrepancy sequence of points?

In the same stream of ideas, given the discrepancy of  the sequence of partitions $\{\pi_n\}$, it would be intresting to estimate the {\it average} speed of the discrepancy of  sequential random reorderings, 

Another interesting problem is to find explicit algorithms to provide low discrepancy sequences of points, given a low discrepancy uniformly distributed sequence of partitions. Some partial results in this direction have been already obtained.

This list of open problems is non exaustive, but it is long enough to confirm Victor Klee's ``conjecture" that the ratio between solved problems and open questions tends to zero, when time goes to infinity.
\bigskip\bigskip\bigskip

Universit\`a della Calabria

Dipartimento di Matematica

87036 Arcavacata di Rende (CS)

Ponte Bucci, cubo 30B

volcic@unical.it
\newpage

\noindent {\bf References}
\bigskip

\noindent [1] ADLER, R.L., FLATTO, L., Uniform distribution of Kakutani's interval splitting procedure, Z. Wahrscheinlichkeitsthorie verw. Gebiete {\bf 38}, (1977) 253-259
\smallskip

\noindent [2] BERGMANN, G., Gleichverteilung von Zerlegunsfolgen nach Kakutani, Diplomarbeit, Erlangen 1984
\smallskip

\noindent [3] BILLINGSLEY, P., Convergence of probability measures, Wiley, 1968
\smallskip

\noindent [4] BRENNAN, M.D., DURRETT, R., Splitting intervals, Annals of Probability, Vol. 14, No. 3 (1986) 1024-1036\smallskip

\noindent [5] CHUNG, K.L. A course in probabilty theory, 2001 Academic Press $3^{rd}$ edition 
\smallskip

\noindent [6] CARBONE, I., VOL\v{C}I\v{C}, A., Kakutani's splitting procedure in higher dimension, Rend. Ist. Matem. Univ. Trieste, Vol. XXXIX (2007), 1-8 
\smallskip

\noindent [7] CARBONE, I., VOL\v{C}I\v{C}, A., A von Neumann theorem for uniformly distributed sequences of partitions, preprint (2007)
\smallskip

\noindent [8] CHERSI, F., VOL\v{C}I\v{C}, A., $\lambda$-equidistributed sequences of partitions and a theorem of the De Bruijn-Post type, Annali Mat. Pura Appl. (IV), Vol. CLXII (1992), 23-32
\smallskip

\noindent [9] DRMOTA, M., TICHY, R.F., Sequences, discrepancies and applications, 1997 Springer LNM {\bf 1651}
\smallskip

\noindent [10] HALTON, J.H. On the efficiency of certain quasi-random sequences of points in evaluating multi-dimensional integrals, Numer. Math., {\bf 2} (1960) 84-90
\smallskip

\noindent [11] HALTON, J.H. Pseudo-random trees: Multiple independent sequence generators for parallel and branching computations, J. Comput. Phys., {\bf 84} (1989) 1-56
\smallskip

\noindent [12] HAMMERSLEY, J.M., Monte Carlo methods for solving multiple problems, Ann. New York Acad. Sci. {\bf 86}, (1960), 844-874
\smallskip

\noindent [13]  HAJIAN, A., ITO, Y., KAKUTANI, S., Invariant measures and orbits of dissipative transformations, Advances in Mathematics {\bf 9}, (1972) 52-65 
\smallskip

\noindent [14] HEWITT, E., SAVAGE, L.J., Symmetric measures on cartesian products, Transactions AMS {\bf 80}, (1955) 470-501
\smallskip

\noindent [15] KAKUTANI, S., A problem on equidistribution on the unit interval [0,1], Measure Theory, Oberwolfach 1975, Springer LNM 541 (1976), 369-375
\smallskip

\noindent [16] KUIPERS, L., NIEDERREITER, H., Uniform distribution of sequences, 1974 Wiley 
\smallskip

\noindent [17] LOOTGIETER, J.C., Sur la r\'epartition des suites de Kakutani (I), Ann. Inst. Henry Poincar\'e, Vol. XIII, No. 4 (1977),  385-410
\smallskip

\noindent [18] LOOTGIETER,  J.C., Sur la r\'epartition des suites de Kakutani (II), Ann. Inst. Henry Poincar\'e, Vol. XIV, No. 3 (1978),  279-302
\smallskip

\noindent [19] PARTHASARATHY, K.R., Probability measures on metric spaces, 1967 Academic Press 
\smallskip

\noindent [20] PREVE, E. Uniformly distributed sequences of partitions, PhD thesis, Milan (2007)
\smallskip

\noindent [21] SCHMIDT, W.M., Irregularities of distribution VII, Acta Arith. {\bf 27}, (1972), 45-50
\smallskip 

\noindent [22] SLUD, E., Entropy and maximal spacing for random partitions, Z. Wahr\-scheinlichkeitsthorie verw. Gebiete {\bf 41}, (1978) 341-352
\smallskip

\noindent [23] SLUD, E., Correction to ``Entropy and maximal spacing for random partitions", Z. Wahrscheinlichkeitsthorie verw. Gebiete {\bf 60}, (1982) 139-141
\smallskip

\noindent [24] VAN DER CORPUT, J.G., Verteilungsfunktionen I, Proc. Akad. Amsterdam {\bf 38,} (1935), 813-821
\smallskip

\noindent [25] VAN ZWET, W.R., A proof of Kakutani's conjecture on random subdivision of longest intervals, Annals of Probability Vol. 6, No. 1 (1978) 133- 137
\smallskip

\noindent [26] 	WEYL, H., \"Uber ein Problem aus dem Gebiete der diophantischen Approximationen, Nach. Ges. Wiss. G\"ottingen, Math.-phys. Kl. (1914), 234-244
\smallskip

\end{document}